\documentclass[12pt]{article}
\usepackage{latexsym,amssymb,amsmath,enumerate,float,geometry,cite}
\newtheorem{theorem}{Theorem}

\newtheorem{lemma}[theorem]{Lemma}

\usepackage{lineno}

\usepackage{setspace}
\onehalfspace
\begin{document}

\title{Exponential Independence}
\author{Simon J\"{a}ger and Dieter Rautenbach}
\date{}
\maketitle
\vspace{-10mm}
\begin{center}
{\small
Institute of Optimization and Operations Research, Ulm University,\\
Ulm, Germany, \texttt{simon.jaeger,dieter.rautenbach@uni-ulm.de}}
\end{center}

\begin{abstract}
For a set $S$ of vertices of a graph $G$, a vertex $u$ in $V(G)\setminus S$, and a vertex $v$ in $S$,
let ${\rm dist}_{(G,S)}(u,v)$ be the distance of $u$ and $v$ in the graph $G-(S\setminus \{ v\})$.
Dankelmann et al. (Domination with exponential decay, Discrete Math. 309 (2009) 5877-5883) 
define $S$ to be an exponential dominating set of $G$
if $w_{(G,S)}(u)\geq 1$ for every vertex $u$ in $V(G)\setminus S$, 
where $w_{(G,S)}(u)=\sum\limits_{v\in S}\left(\frac{1}{2}\right)^{{\rm dist}_{(G,S)}(u,v)-1}$.
Inspired by this notion, 
we define $S$ to be an exponential independent set of $G$
if $w_{(G,S\setminus \{ u\})}(u)<1$ for every vertex $u$ in $S$,
and the exponential independence number $\alpha_e(G)$ of $G$ as the maximum order of an exponential independent set of $G$.

Similarly as for exponential domination, the non-local nature of exponential independence leads to many interesting effects and challenges.
Our results comprise exact values for special graphs as well as tight bounds and the corresponding extremal graphs.
Furthermore, we characterize all graphs $G$ for which $\alpha_e(H)$ equals the independence number $\alpha(H)$ for every induced subgraph $H$ of $G$,
and we give an explicit characterization of all trees $T$ with $\alpha_e(T)=\alpha(T)$.
\end{abstract}

\pagebreak

\section{Introduction}

Independence in graphs is one of the most fundamental and well-studied concepts in graph theory.
In the present paper we propose and study a version of independence 
where the influence of vertices decays exponentially with respect to distance.
This new notion is inspired by the {\it exponential domination number},
which was introduced by Dankelmann et al. \cite{ddems} and recently studied in \cite{abcvy,aa,bor1,bor2}.
Somewhat related parameters are the 
well-known {\it (distance) packing numbers} \cite{jr,hen,mm}
and the {\it influence numbers} \cite{dll,hed}.

We consider finite, simple, and undirected graphs, and use standard terminology.
The vertex set and the edge set of a graph $G$ are denoted by $V(G)$ and $E(G)$, respectively.
The order $n(G)$ of $G$ is the number of vertices of $G$.
The distance ${\rm dist}_G(u,v)$ between two vertices $u$ and $v$ in a graph $G$ 
is the minimum number of edges of a path in $G$ between $u$ and $v$.
If no such path exists, then let ${\rm dist}_G(u,v)=\infty$.
The diameter ${\rm diam}(G)$ of $G$ is the maximum distance between vertices of $G$.
A set of pairwise non-adjacent vertices of $G$ is an independent set of $G$,
and the maximum order of an independent set of $G$ is the independence number $\alpha(G)$ of $G$.

Let $S$ be a set of vertices of $G$.
For two vertices $u$ and $v$ of $G$,
let ${\rm dist}_{(G,S)}(u,v)$ be the minimum number of edges of a path $P$ in $G$ between $u$ and $v$ 
such that $S$ contains exactly one endvertex of $P$ but no internal vertex of $P$.
If no such path exists, then let ${\rm dist}_{(G,S)}(u,v)=\infty$.
Note that, if $u$ and $v$ are distinct vertices in $S$, 
then ${\rm dist}_{(G,S)}(u,u)=0$ and ${\rm dist}_{(G,S)}(u,v)=\infty$.
For a vertex $u$ of $G$, let 
\begin{eqnarray}\label{ew}
w_{(G,S)}(u)=\sum\limits_{v\in S}\left(\frac{1}{2}\right)^{{\rm dist}_{(G,S)}(u,v)-1},
\end{eqnarray}
where $\left(\frac{1}{2}\right)^{\infty}=0$.
Note that $w_{(G,S)}(u)=2$ for $u\in S$.

Dankelmann et al. \cite{ddems} define a set $S$ of vertices to be {\it exponential dominating} if
$$\mbox{$w_{(G,S)}(u)\geq 1$ for every vertex $u$ in $V(G)\setminus S$,}$$ 
and the {\it exponential domination number $\gamma_e(G)$} of $G$ 
as the minimum order of an exponential dominating set.
Analogously, we define $S$ to be {\it exponential independent} if 
$$\mbox{$w_{(G,S\setminus \{ u\})}(u)<1$ for every vertex $u$ in $S$,}$$ 
that is, the accumulated exponentially decaying influence $w_{(G,S\setminus \{ u\})}(u)$ of the remaining vertices in $S\setminus \{ u\}$ 
that arrives at any vertex $u$ in $S$ is strictly less than $1$.
Let the {\it exponential independence number $\alpha_e(G)$} of $G$ 
be the maximum order of an exponential independent set.
An (exponential) independent set of maximum order is {\it maximum}.

Our results comprise exact values for special graphs as well as tight bounds and the corresponding extremal graphs.
Furthermore, we characterize all graphs $G$ for which $\alpha_e(H)$ equals the independence number $\alpha(H)$ for every induced subgraph $H$ of $G$,
and we give an explicit characterization of all trees $T$ with $\alpha_e(T)=\alpha(T)$.
We conclude with several open problems.

\section{Results}

We start with some elementary observations concerning exponential independence.
Clearly, every exponential independent set is independent, 
which immediately implies (i) of the following theorem.
The quantity $w_{(G,S\setminus \{ u\})}(u)$ does not behave monotonously with respect to the removal of vertices from $S$.
Indeed, if $G$ is a star $K_{1,n-1}$ with center $v$, and $S=V(G)$ for instance,
then $w_{(G,S\setminus \{ u\})}(u)=1$ for every endvertex $u$ of $G$
but $w_{(G,S\setminus \{ u,v\})}(u)=\frac{n-2}{2}$,
which can be smaller or bigger than $1$.
In view of this observation part (iii) of the following theorem is slightly surprising.

\begin{theorem}\label{theorem1}
Let $G$ be a graph.
\begin{enumerate}[(i)]
\item $\alpha_e(G)\leq \alpha(G)$.
\item If $H$ is a subgraph of $G$ and $S\subseteq V(H)$ is an exponential independent set of $G$,
then $S$ is an exponential independent set of $H$.
\item A subset of an exponential independent set of $G$ is an exponential independent set of $G$.
\end{enumerate}
\end{theorem}
{\it Proof:} (i) follows from the above observation.
Since ${\rm dist}_{(G,S\setminus \{ u\})}(u,v)\leq {\rm dist}_{(H,S\setminus \{ u\})}(u,v)$ for every two vertices $u$ and $v$ in $S$, 
(ii) follows immediately from (\ref{ew}).
We proceed to the proof of (iii).
Let $S$ be an exponential independent set of $G$.
Let $u$ and $v$ be distinct vertices in $S$.
In order to complete the proof, 
it suffices to show 
\begin{eqnarray}\label{e1}
w_{(G,S\setminus \{ u,v\})}(u)\leq w_{(G,S\setminus \{ u\})}(u).
\end{eqnarray}
For 
\begin{eqnarray*}
S_\infty & =& \{ w\in S\setminus \{ u,v\}:{\rm dist}_{(G,S\setminus \{ u,v\})}(u,w)=\infty\},\\
S_=&=& \{ w\in S\setminus \{ u,v\}:{\rm dist}_{(G,S\setminus \{ u,v\})}(u,w)={\rm dist}_{(G,S\setminus \{ u\})}(u,w)<\infty\},\mbox{ and}\\
S_>&=& \{ w\in S\setminus \{ u,v\}:{\rm dist}_{(G,S\setminus \{ u,v\})}(u,w)<{\rm dist}_{(G,S\setminus \{ u\})}(u,w)\},
\end{eqnarray*}
we have $S=\{ u,v\}\cup S_=\cup S_>\cup S_\infty$.
If $S_>=\emptyset$, then (\ref{e1}) follows immediately from (\ref{ew}).
Hence, we may assume that $S_>\not=\emptyset$.
Let $T$ be a subtree of $G$ rooted in $u$ such that 
\begin{itemize}
\item $S_=\cup S_>$ is the set of all leaves of $T$, 
\item ${\rm dist}_T(u,w)={\rm dist}_{(G,S\setminus \{ u,v\})}(u,w)$ for every $w\in S_=\cup S_>$, and
\item $v$ is not an ancestor within $T$ of any vertex in $S_=$.
\end{itemize}
Such a tree can easily be extracted from the union of paths $P_w$ for $w\in S_=\cup S_>$,
where $P_w$ is a path of length ${\rm dist}_{(G,S\setminus \{ u,v\})}(u,w)$
between $w$ and $u$ that intersects $S\setminus \{ u,v\}$ only in $w$,
and that avoids $v$ if $w\in S_=$.
Since $S_>\not=\emptyset$, the vertex $v$ belongs to $T$, 
and the set of leaves of $T$ that are descendants of $v$ is exactly $S_>$.
The conditions imposed on $T$ easily imply ${\rm dist}_T(u,v)={\rm dist}_{(G,S\setminus \{ u\})}(u,v)$.
Let $T_>$ be the subtree of $T$ rooted in $v$ that contains $v$ and all its descendants within $T$.
Since $S$ is exponential independent, we obtain $w_{(T_>,S_>)}(v)\leq w_{(G,S\setminus \{ v\})}(v)<1$, which implies
\begin{eqnarray*}
w_{(G,S\setminus \{ u,v\})}(u) 
& = & w_{(T,S_=)}(u)+w_{(T,S_>)}(u)\\
& = & w_{(T,S_=)}(u)+\left(\frac{1}{2}\right)^{{\rm dist}_T(u,v)}w_{(T_>,S_>)}(v)\\
& < & w_{(T,S_=)}(u)+\left(\frac{1}{2}\right)^{{\rm dist}_{(G,S\setminus \{ u\})}(u,v)}\\
& = & \sum_{w\in S_=}\left(\frac{1}{2}\right)^{{\rm dist}_{(G,S\setminus \{ u\})}(u,w)-1}+\left(\frac{1}{2}\right)^{{\rm dist}_{(G,S\setminus \{ u\})}(u,v)}\\
& <  & \sum_{w\in S_=\cup \{ v\}}\left(\frac{1}{2}\right)^{{\rm dist}_{(G,S\setminus \{ u\})}(u,w)-1}\\
& \leq & w_{(G,S\setminus \{ u\})}(u),
\end{eqnarray*}
which completes the proof. $\Box$

\medskip

\noindent Our next result is a lower bound on the exponential independence number,
for which we are able to characterize all extremal trees.

\begin{theorem}\label{theorem2}
If $G$ is a connected graph of order $n$ and diameter ${\rm diam}$, then 
\begin{eqnarray}\label{elower}
\alpha_e(G)\geq \frac{2{\rm diam}+2}{5}.
\end{eqnarray}
Furthermore, if $G$ is a tree, 
then (\ref{elower}) holds with equality 
if and only if $G$ is a path and $n$ is a multiple of $5$. 
\end{theorem}
{\it Proof:} Let $P:v_0v_1\ldots v_{\rm diam}$ be a shortest path of length ${\rm diam}$ in $G$.
Let 
$$S=\left\{v_{5i}:  i\in\left\{0,\ldots, \left \lfloor\frac{{\rm diam}}{5}\right\rfloor\right\}\right\}
\cup\left\{v_{5i+2}: i\in\left\{0,\ldots, \left \lfloor\frac{{\rm diam}-2}{5}\right\rfloor\right\}\right\}.$$
Let $v_i\in S$.
Since $P$ is a shortest path, we have
${\rm dist}_{(G,S\setminus \{ v_i\})}(v_i,v_j)\geq |j-i|$ for every $v_j$ in $S\setminus \{ v_i\}$.
By construction, the set $S$ contains no neighbor of $v_i$,
and
$S$ contains at most one of the two vertices $v_{i-k}$ and $v_{i+k}$ for every integer $k$ at least $2$.
This implies 
$w_{(G,S\setminus \{ v_i\})}(v_i)<\sum\limits_{k=2}^{\infty}\left(\frac{1}{2}\right)^{k-1}=1$.
Hence, $S$ is an exponential independent set of $G$, and
\begin{eqnarray*}
\alpha_e(G) & \geq & |S|
= 1+\left\lfloor\frac{{\rm diam}}{5}\right\rfloor+1+\left\lfloor\frac{{\rm diam}-2}{5}\right\rfloor
\geq \frac{2{\rm diam}+2}{5}.
\end{eqnarray*}
Now, let $G$ be a path and let $n$ be a multiple of $5$, that is, $G=P_n$.
It is easy to verify that $\alpha_e(P_5)=2=\frac{2{\rm diam}+2}{5}$.
Furthermore, if $n>5$ and $S$ is a maximum exponential independent set of $G$,
then $S\cap \{ v_0,v_1,v_2,v_3,v_4\}$ is an exponential independent set of $P_5$
and $S\setminus \{ v_0,v_1,v_2,v_3,v_4\}$ is an exponential independent set of $P_{n-5}$.
By an inductive argument, we obtain,
$$\frac{2n}{5}\leq \alpha_e(G)=\alpha_e(P_n)\leq \alpha_e(P_5)+\alpha_e(P_{n-5})=2+\frac{2(n-5)}{5}=\frac{2n}{5},$$
which implies that paths whose order is a multiple of $5$ satisfy (\ref{elower}) with equality.

Finally, let $G$ be a tree with $\alpha_e(G)=\frac{2{\rm diam}+2}{5}$, and let $P$ be as above.
Since $\frac{2{\rm diam}+2}{5}$ is an integer, the order ${\rm diam}+1$ of $P$ is a multiple of $5$.
Suppose that $G$ is distinct from $P$.
This implies that there is some vertex $v_k$ of $P$ that has a neighbor $u$ that does not belong to $P$.
Let $k=5r+s$ for some $s\in \{ 0,1,2,3,4\}$.
By symmetry, we may assume that $s\in \{ 0,1,2\}$.
Let 
\begin{eqnarray*}
S_0 & = & \{ v_i:i\in \{ 0,\ldots,5r-1\}\mbox{ with $i$ mod $5$}\in \{ 0,2\}\}\\
&& \cup \{ v_i:i\in \{ 5r+5,\ldots,{\rm diam}+1\}\mbox{ with $i$ mod $5$}\in \{ 2,4\}\}
\end{eqnarray*}
If $s=0$, then let $S=\{ v_{k+1},v_{k+4},u\}\cup S_0$, and
if $s\in \{ 1,2\}$, then let $S=\{ v_k,v_{k+4},u\}\cup S_0$.
The set $S$ is an exponential independent set of $G$ of order more than $\frac{2{\rm diam}+2}{5}$,
which is a contradiction.
Hence, $G$ is a path and $n$ is a multiple of $5$,
which completes the proof. $\Box$

\medskip

\noindent For later reference, we include a fundamental lemma from \cite{bor2}.
Recall that a full binary tree is a rooted tree in which each vertex has either no or exactly two children.

\begin{lemma}[Bessy et al. \cite{bor2}]\label{lemma1}
Let $G$ be a graph of maximum degree at most $3$,
and let $S$ be a set of vertices of $G$.

If $u$ is a vertex of degree at most $2$ in $G$, then $w_{(G,S)}(u)\leq 2$ with equality if and only if 
$u$ is contained in a subgraph $T$ of $G$ that is a tree, 
such that rooting $T$ in $u$ yields a full binary tree and $S\cap V(T)$ is exactly the set of leaves of $T$.
\end{lemma}
Our next result concerns the exponential independence numbers of some special graphs.

\begin{theorem}\label{theorem3}
\begin{enumerate}[(i)]
\item If $P_n$ is the path of order $n$, then $\alpha_e(P_n)=\left\lceil\frac{2n}{5}\right\rceil$.
\item If $C_n$ is the path of order $n$ at least $5$, then $\alpha_e(C_n)=\left\lfloor\frac{2n}{5}\right\rfloor$.
\item If $T$ is a full binary tree of order $n$, then $\alpha_e(T)=\frac{n+1}{2}$.
Furthermore, the set of leaves of $T$ is the unique maximum exponential independent set of $T$.
\end{enumerate}
\end{theorem}
{\it Proof:} 
(i) By Theorem \ref{theorem2}, 
$\alpha_e(P_n)$ is at least $\left\lceil\frac{2n}{5}\right\rceil$.
For $n\leq 5$, it is easy to verify that $\alpha_e(P_n)$ is also at most $\left\lceil\frac{2n}{5}\right\rceil$.
Now, let $n>5$. Let $P_n$ be the path $v_0v_1\ldots v_{n-1}$.
Let $S$ be a maximum exponential independent set of $P_n$.
Since $S\cap \{ v_0,v_1,v_2,v_3,v_4\}$ is an exponential independent set of $P_5$
and $S\setminus \{ v_0,v_1,v_2,v_3,v_4\}$ is an exponential independent set of $P_{n-5}$,
we obtain, by an inductive argument,
$$\alpha_e(P_n)\leq \alpha_e(P_5)+\alpha_e(P_{n-5})=2+\left\lceil\frac{2(n-5)}{5}\right\rceil=\left\lceil\frac{2n}{5}\right\rceil.$$
(ii) For $5\leq n\leq 9$, it is easy to verify that $\alpha_e(C_n)=\left\lfloor\frac{2n}{5}\right\rfloor$.
For $n\geq 10$, a similar argument as for the paths implies $\alpha_e(C_n)=2+\alpha_e(C_{n-5})$, 
and an inductive argument yields (ii).

\medskip

\noindent (iii) Clearly, we may assume $n>3$.
Let $L$ be the set of leaves of $T$.
Note that $n=2|L|-1$.
Let $v\in L$ and let $u$ be the parent of $u$ in $T$.
If $w_{(T,L\setminus \{ v\})}(v)\geq 1$, then $w_{(T,L\setminus \{ v\})}(u)\geq 2$,
and Lemma \ref{lemma1} implies that rooting $T-v$ in $u$ yields a full binary tree.
This implies the contradiction that $T$ only has vertices of degree $1$ and $3$,
while the root of $T$ has degree $2$.
Hence, $L$ is an exponential independent set of $T$, which implies $\alpha_e(T)\geq |L|=\frac{n+1}{2}$.

Suppose that $T$ is a full binary tree of minimum order $n$ such that
either $\alpha_e(T)>\frac{n+1}{2}$
or $\alpha_e(T)=\frac{n+1}{2}$ but $T$ has a maximum exponential independent set distinct from $L$.
In both cases, $T$ has a maximum exponential independent set $S$ with $S\setminus L\not=\emptyset$.
Let $v$ be a vertex in $S\setminus L$ at maximum distance from the root of $T$.
Let $w$ and $w'$ be the two children of $v$ in $T$.
Since $w_{(T,S\setminus \{ v\})}(v)<1$,
Lemma \ref{lemma1} implies that $L\setminus S$ contains 
at least one leaf of $T$ that is either $w$ or a descendant of $w$
as well as 
at least one leaf of $T$ that is either $w'$ or a descendant of $w'$.
Hence, if $\ell_v$ leaves of $T$ are descendants of $v$, then $S$ contains at most $\ell_v-2$ descendants of $v$.
Let $T'$ arise from $T$ by removing all descendants of $v$.
Since $T'$ is a full binary tree of smaller order than $T$,
the choice of $T$ implies that $\alpha_e(T')=\frac{n(T')+1}{2}$.
Note that $S\cap V(T')$ is an exponential independent set of $T'$, 
and that $v$ has exactly $2\ell_v-2$ descendants.
Therefore,
\begin{eqnarray*}
|S| & \leq & |S\setminus V(T')|+|S\cap V(T')|\\
& \leq & (\ell_v-2)+|S\cap V(T')|\\
& \leq & (\ell_v-2)+\frac{n(T')+1}{2}\\
&= & (\ell_v-2)+\frac{n-(2\ell_v-2)+1}{2}\\
& < & \frac{n+1}{2},
\end{eqnarray*}
which is a contradiction.
$\Box$

\medskip

\noindent Our next result is an upper bound on the exponential independence number,
for which we achieve a full characterization of the extremal graphs.

\begin{theorem}\label{theorem2b}
If $G$ is a connected graph of order $n$, then 
$$
\alpha_e(G)\leq \frac{n+1}{2}
$$
with equality if and only if $G$ is a full binary tree.
\end{theorem}
{\it Proof:} 
We show the upper bound by induction on $n$.
By Theorem \ref{theorem1}(ii), we may assume that $G$ is a tree $T$.
If $n=1$, then $\alpha_e(T)=1=\frac{n+1}{2}$.
Now, let $n\geq 2$, and let $S$ be a maximum exponential independent set of $T$.
We root $T$ in some vertex $r$.
Let $v$ be a vertex in $S$ at maximum distance from $r$.
If $v=r$, then $|S|=1$, and the statement holds.
Hence, we may assume that $v$ and $r$ are distinct.
Let $u$ be the parent of $v$.

First, we assume that $v$ is the only descendant of $u$ that belongs to $S$.
Let $T'$ arise from $T$ by removing $u$ together with all descendants of $u$, and let $S'=S\setminus \{ v\}$.
Clearly, $S'$ is an exponential independent set of the tree $T'$, and we obtain, by induction,
$$\alpha_e(T)
=|S|
=|S'|+1
\leq \alpha_e(T')+1
\leq \frac{n(T')+1}{2}+1
\leq \frac{n(T)+1}{2}.$$
Next, we assume that $S$ contains some descendant of $u$ distinct from $v$.
Let $S_u$ be the set of descendants of $u$ that belong to $S$.
By the choice of $v$, all vertices in $S_u$ are children of $u$.
Since $S$ is exponential independent, we obtain $|S_u|=2$, and $u\not\in S$.
Let $T''$ arise from $T$ by removing all descendants of $u$, and let $S''=(S\setminus S_u)\cup \{ u\}$.
If $w_{(T'',S''\setminus \{ u\})}(u)\geq 1$, 
then $w_{(T,S\setminus \{ v\})}(v)\geq \frac{1}{2}w_{(T'',S''\setminus \{ u\})}(u)+\frac{1}{2}\geq 1$,
which is a contradiction.
If $w_{(T'',S''\setminus \{ w\})}(w)\geq 1$ for some $w\in S''\setminus \{ u\}$, 
then ${\rm dist}_{(T'',S''\setminus \{ w\})}(w,u)={\rm dist}_{(T,S\setminus \{ w\})}(w,x)-1$ 
for every $x\in S_u$ and $|S_u|=2$ imply
\begin{eqnarray*}
w_{(T,S\setminus \{ w\})}(w)
&=& w_{(T'',S''\setminus \{ w\})}(w)
-\left(\frac{1}{2}\right)^{{\rm dist}_{(T'',S''\setminus \{ w\})}(w,u)-1} 
+\sum_{x\in S_u}\left(\frac{1}{2}\right)^{{\rm dist}_{(T,S\setminus \{ w\})}(w,x)-1}\\ 
&=& w_{(T'',S''\setminus \{ w\})}(w)\\
&\geq & 1,
\end{eqnarray*}
which is a contradiction.
Hence, $S''$ is an exponential independent set of $T''$, and we obtain, by induction,
$$\alpha_e(T)
=|S|
=|S''|+1
\leq \alpha_e(T'')+1
\leq \frac{n(T'')+1}{2}+1
\leq \frac{n(T)+1}{2},$$
which completes the proof of the upper bound.

Next, we show that we have equality if and only if $G$ is a full binary tree.
By Theorem \ref{theorem3}(iii),
we only need to show that every connected graph $G$ with $\alpha_e(G)=\frac{n+1}{2}$ is a full binary tree.
Therefore, suppose that $G$ is a connected graph of minimum order $n$ with $\alpha_e(G)=\frac{n+1}{2}$ that is not a full binary tree.

Let $T$ be a spanning tree of $G$.
We will show first that $T$ is a full binary tree.
By Theorem \ref{theorem1}(ii), we have $\frac{n+1}{2}=\alpha_e(G)\leq \alpha_e(T)\leq \frac{n+1}{2}$,
which implies $\alpha_e(T)=\frac{n+1}{2}$.
Let $S$ be a maximum exponential independent set of $G$, and, hence, also of $T$.
If the diameter of $T$ is at most $2$, then it is easy to see that either $\alpha_e(G)\neq\frac{n+1}{2}$ or $G$ is a full binary tree,
that is, the diameter of $T$ is at least $3$.
Let $w$ be the endvertex of a longest path $P$ in $T$.
Let $v$ be the neighbor of $w$, and let $u$ be the neighbor of $v$ on $P$ that is distinct from $w$.
Let $T'=T-(N_T[v]\setminus \{ u\})$, and let $S'=S\cap V(T')$.
Note that all vertices in $N_T(v)\setminus \{ u\}$ are endvertices of $T$.

First, we assume that $v$ has degree $2$ in $T$.
Note that $S'$ is an exponential independent set of $T'$, 
the set $S$ contains at most one of the two vertices $v$ and $w$, and $n(T')=n-2$.
This implies 
$\frac{n+1}{2}=\alpha_e(T)=|S|\leq |S'|+1\leq \alpha_e(T')+1\leq \frac{n(T')+1}{2}+1=\frac{n+1}{2}$,
which implies that $\alpha_e(T')=|S'|=\frac{n(T')+1}{2}$, and that $S$ contains either $w$ or $v$.
By the choice of $G$, this implies that $T'$ is a full binary tree.
By Theorem \ref{theorem3}(iii), the set $S'$ is exactly the set of leaves of $T'$.
If $u$ is the root of $T'$, then $T$ is a full binary tree with root $v$, which is a contradiction.
Hence, $u$ is not the root of $T'$.
Let $u'$ be a leaf of $T'$ that is either $u$ or a descendant of $u$ in $T'$.
Let $u'$ have distance $d$ from the root of $T'$.
Note that the distance between $w$ and $u'$ is at most $d+1$.
Now, Lemma \ref{lemma1} implies the contradiction 
$w_{(T,S\setminus \{ u'\})}(u')\geq \left(\frac{1}{2}\right)^{(d+1)-1}+\sum\limits_{i=1}^d\left(\frac{1}{2}\right)^i=1$.
Hence, $v$ has degree at least $3$ in $T$.

If $S$ contains at most one vertex from $N_T[v]\setminus \{ u\}$, then $n(T')\leq n-3$ implies the contradiction
$\frac{n+1}{2}=\alpha_e(T)=|S|\leq |S'|+1\leq \alpha_e(T')+1\leq \frac{n(T')+1}{2}+1=\frac{n}{2}$.
It follows easily that $S$ contains exactly two vertices from $N_T(v)\setminus \{ u\}$ but not $v$.
Let $T''=T-(N_T(v)\setminus \{ u\})$, and let $S''=S'\cup \{ v\}$.
Arguing as before, it follows that $S''$ is an exponential independent set of $T''$.
Since $n(T'')\leq n-2$, we obtain
$\frac{n+1}{2}=\alpha_e(T)=|S|=|S''|+1\leq \alpha_e(T'')+1\leq \frac{n(T'')+1}{2}+1\leq\frac{n+1}{2}$,
which implies $\alpha_e(T'')=\frac{n(T'')+1}{2}$ and $n(T'')=n-2$.
By the choice of $G$, it follows that $T''$ is a full binary tree, and that $S''$ is a maximum exponential independent set of $T''$.
Since $v\in S''$, Theorem \ref{theorem3} implies that $v$ is a leaf of $T''$.
Now, also in this case, the tree $T$ is a full binary tree.

Since $T$ was an arbitrary spanning tree of $G$, it follows that every spanning tree of $G$ is a full binary tree.
This easily implies that $G=T$, that is, $G$ is a full binary tree, which completes the proof. $\Box$

\medskip

\noindent Theorem \ref{theorem2} implies that 
$\alpha_e(G)$ is at least $\Omega(\log_2(n(G)))$ for every connected cubic graph $G$.
We conjecture that $\alpha_e(G)$ actually grow much faster than $\log_2(n(G))$.
At least for subcubic trees, we obtain the following linear lower bound.

\begin{theorem}\label{theorem4}
If $T$ is a tree of order $n$ and maximum degree at most $3$, then $\alpha_e(T)\geq \frac{2n+8}{13}$.
\end{theorem}
{\it Proof:} Clearly, we may assume that $n>3$.
Let $T$ have $n_i$ vertices of degree $i$ for $i\in [3]$.
Note that $n_1\geq n_3+2$.

If $n_2>0$, then let $S_1$ be the set of all leaves of $T$, and, 
if $n_2=0$, then let $S_1$ be the set of all leaves of $T$ except for exactly one.
Arguing as in the proof of Theorem \ref{theorem3}(iii), it follows that $S_1$ is an exponential independent set in $T$,
which implies $\alpha_e(T)\geq n_1-1\geq n_3+1$.

Let $V_3$ be the set of vertices of degree $3$, and let $T'=T-N_T[V_3]$.
Note that $T'$ is a union of paths, and that $n(T')\geq n-4n_3$.
By Theorem \ref{theorem3}(i), 
the forest $T'$ has an exponential independent set $S_2$ of order at least $\frac{2n(T')}{5}\geq \frac{2n-8n_3}{5}$.
We will show that $S_2$ is also exponential independent within $T$.
Therefore, let $u$ be a vertex of degree $1$ in $T'$ that has a neighbor $v$ in $V(T)\setminus V(T')$.
By construction, $u$ and $v$ have degree $2$ in $T$, and $v$ has a neighbor $w$ of degree $3$ in $T$.
Let $T_w$ be the component of $T-v$ that contains $w$, and let $S_w=S_2\cap V(T_w)$.
If $w_{(T,S_w)}(u)\geq \frac{1}{2}$, then $w_{(T_w,S_w)}(w)\geq 2$.
By Lemma \ref{lemma1}, this implies that $S_w$, and hence $S_2$, intersects $N_T[V_3]$, which is a contradiction.
Hence, $w_{(T,S_w)}(u)<\frac{1}{2}$.
Similarly, if $u$ is a vertex of degree $0$ in $T'$, 
then $w_{(T,S_2\setminus \{ u\})}(u)<\frac{1}{2}$ if $u$ has degree $1$ in $T$,
and 
$w_{(T,S_2\setminus \{ u\})}(u)<\frac{1}{2}+\frac{1}{2}$ if $u$ has degree $2$ in $T$.
If $P=v_0\ldots v_\ell$ is a component of $T'$ with $|V(P)\cap S_2|\geq 2$, 
and $v_i\in S_2$ is such that $S\cap \{ v_0,\ldots,v_{i-1}\}=\emptyset$,
then $w_{(T',S_2\setminus \{ v_i\})}(v_i)\leq \frac{1}{2}$.
Combining these observations, it follows easily that $S_2$ is an exponential independent set in $T$, 
which implies $\alpha_e(T)\geq \frac{2n-8n_3}{5}$.

Altogether, we obtain
$\alpha_e(T)\geq \max\left\{ n_3+1,\frac{2n-8n_3}{5}\right\}\geq \frac{2n+8}{13},$
which completes the proof. $\Box$

\medskip

\noindent After the above bounds, exact values, and extremal graphs, we consider graphs $G$ with $\alpha_e(G)=\alpha(G)$.
We achieve full characterizations of all graphs for which every induced subgraph has this property,
and also of all trees that have this property.

Recall that the {\it bull} is the unique graph $B$ of order $5$ with degree sequence $1,1,2,3,3$.

\begin{theorem}\label{theorem5}
If $G$ is a graph, then $\alpha_e(H)=\alpha(H)$ for every induced subgraph $H$ of $G$ if and only if $G$ is $\{ K_{1,3},P_5,B\}$-free.
\end{theorem}
{\it Proof:} If $H\in \{ K_{1,3},P_5,B\}$, then $\alpha_e(H)=2<3=\alpha(H)$, which implies the necessity.
In order to show the sufficiency, let $G$ be a $\{ K_{1,3},P_5,B\}$-free graph.
It suffices to show that $\alpha_e(G)=\alpha(G)$.
Let $S$ be a maximum independent set of $G$.
If $|S|\leq 2$, then $S$ is also exponential independent, which implies $\alpha_e(G)=\alpha(G)$.
Hence, we may assume that $|S|\geq 3$.
Possibly iteratively replacing elements of $S$ by one of their neighbors, 
we may assume that $S$ contains two vertices $u$ and $v$ at distance $2$.
Suppose that $S\setminus \{ u,v\}$ contains a vertex $w$ at distance $2$ from $u$.
If $u$, $v$, and $w$ have a common neighbor, 
then the independence of $S$ implies that $G$ contains $K_{1,3}$ as an induced subgraph,
which is a contradiction. 
Therefore, if $uv'v$ and $uw'w$ are shortest paths in $G$, 
then $v'\not=w'$ and $vw',wv'\not\in E(G)$,
which implies the contradiction that $\{ u,v,w,v',w'\}$ induces $P_5$ or $B$.
Hence, we may assume, by symmetry, that for no vertex $x$ in $S$,
there are two vertices in $S$ at distance $2$ from $x$.
Let $w\in S\setminus \{ u,v\}$. 
Since $G$ is $P_5$-free, the distance of $u$ and $w$ is $3$.
Let $uv'v$ and $uw'w''w$ be shortest paths in $G$.
Note that $v$ is not adjacent to $w''$.
If $v'=w'$, then $\{ u,v,v',w''\}$ induces $K_{1,3}$, which is a contradiction.
Hence, $v'\not=w'$.
By symmetry, we may assume that $v$ is not adjacent to $w'$ and that $v'$ is not adjacent to $w''$.
Now, $\{ u,v,v',w',w''\}$ induces $P_5$ or $B$, 
which is a contradiction and completes the proof.
$\Box$

\medskip

\noindent We proceed to the trees $T$ with $\alpha_e(T)=\alpha(T)$.

\begin{figure}[H]
\begin{center}
{\footnotesize
%TeXCAD Picture [1.pic]. Options:
%\grade{\on}
%\emlines{\off}
%\epic{\off}
%\beziermacro{\on}
%\reduce{\on}
%\snapping{\on}
%\pvinsert{% Your \input, \def, etc. here}
%\quality{8.000}
%\graddiff{0.005}
%\snapasp{1}
%\zoom{5.6569}
\unitlength 0.8mm % = 2.845pt
\linethickness{0.4pt}
\ifx\plotpoint\undefined\newsavebox{\plotpoint}\fi % GNUPLOT compatibility
\begin{picture}(186,25)(0,0)
\put(5,5){\circle*{2}}
\put(155,5){\circle*{2}}
\put(105,5){\circle*{2}}
\put(65,5){\circle*{2}}
\put(35,5){\circle*{2}}
\put(145,5){\circle*{2}}
\put(95,5){\circle*{2}}
\put(55,5){\circle*{2}}
\put(25,5){\circle*{2}}
\put(165,5){\circle*{2}}
\put(115,5){\circle*{2}}
\put(75,5){\circle*{2}}
\put(175,5){\circle*{2}}
\put(125,5){\circle*{2}}
\put(185,5){\circle*{2}}
\put(5,15){\circle*{2}}
\put(155,15){\circle*{2}}
\put(105,15){\circle*{2}}
\put(65,15){\circle*{2}}
\put(35,15){\circle*{2}}
\put(145,15){\circle*{2}}
\put(95,15){\circle*{2}}
\put(55,15){\circle*{2}}
\put(25,15){\circle*{2}}
\put(165,15){\circle*{2}}
\put(115,15){\circle*{2}}
\put(75,15){\circle*{2}}
\put(175,15){\circle*{2}}
\put(125,15){\circle*{2}}
\put(185,15){\circle*{2}}
\put(5,5){\line(0,1){10}}
\put(155,5){\line(0,1){10}}
\put(105,5){\line(0,1){10}}
\put(65,5){\line(0,1){10}}
\put(35,5){\line(0,1){10}}
\put(145,5){\line(0,1){10}}
\put(95,5){\line(0,1){10}}
\put(55,5){\line(0,1){10}}
\put(25,5){\line(0,1){10}}
\put(165,5){\line(0,1){10}}
\put(115,5){\line(0,1){10}}
\put(75,5){\line(0,1){10}}
\put(175,5){\line(0,1){10}}
\put(125,5){\line(0,1){10}}
\put(185,5){\line(0,1){10}}
\put(145,15){\line(1,0){10}}
\put(95,15){\line(1,0){10}}
\put(55,15){\line(1,0){10}}
\put(25,15){\line(1,0){10}}
\put(155,15){\line(1,0){10}}
\put(105,15){\line(1,0){10}}
\put(65,15){\line(1,0){10}}
\put(165,15){\line(1,0){10}}
\put(115,15){\line(1,0){10}}
\put(175,15){\line(1,0){10}}
\put(5,19){\makebox(0,0)[cc]{$x_1$}}
\put(155,19){\makebox(0,0)[cc]{$x_2$}}
\put(105,19){\makebox(0,0)[cc]{$x_2$}}
\put(65,19){\makebox(0,0)[cc]{$x_2$}}
\put(35,19){\makebox(0,0)[cc]{$x_2$}}
\put(145,19){\makebox(0,0)[cc]{$x_1$}}
\put(95,19){\makebox(0,0)[cc]{$x_1$}}
\put(55,19){\makebox(0,0)[cc]{$x_1$}}
\put(25,19){\makebox(0,0)[cc]{$x_1$}}
\put(165,19){\makebox(0,0)[cc]{$x_3$}}
\put(115,19){\makebox(0,0)[cc]{$x_3$}}
\put(75,19){\makebox(0,0)[cc]{$x_3$}}
\put(175,19){\makebox(0,0)[cc]{$x_4$}}
\put(125,19){\makebox(0,0)[cc]{$x_4$}}
\put(185,19){\makebox(0,0)[cc]{$x_5$}}
\put(5,1){\makebox(0,0)[cc]{$y_1$}}
\put(155,1){\makebox(0,0)[cc]{$y_2$}}
\put(105,1){\makebox(0,0)[cc]{$y_2$}}
\put(65,1){\makebox(0,0)[cc]{$y_2$}}
\put(35,1){\makebox(0,0)[cc]{$y_2$}}
\put(145,1){\makebox(0,0)[cc]{$y_1$}}
\put(95,1){\makebox(0,0)[cc]{$y_1$}}
\put(55,1){\makebox(0,0)[cc]{$y_1$}}
\put(25,1){\makebox(0,0)[cc]{$y_1$}}
\put(165,1){\makebox(0,0)[cc]{$y_3$}}
\put(115,1){\makebox(0,0)[cc]{$y_3$}}
\put(75,1){\makebox(0,0)[cc]{$y_3$}}
\put(175,1){\makebox(0,0)[cc]{$y_4$}}
\put(125,1){\makebox(0,0)[cc]{$y_4$}}
\put(185,1){\makebox(0,0)[cc]{$y_5$}}
\put(165,25){\makebox(0,0)[cc]{$T_1(5)$}}
\put(110,25){\makebox(0,0)[cc]{$T_1(4)$}}
\put(65,25){\makebox(0,0)[cc]{$T_1(3)$}}
\put(30,25){\makebox(0,0)[cc]{$T_1(2)$}}
\put(5,25){\makebox(0,0)[cc]{$T_1(1)$}}
\end{picture}
}
\end{center}
\caption{The trees $T_1(k)$ for $k\in [5]$.}\label{fig1}
\end{figure}
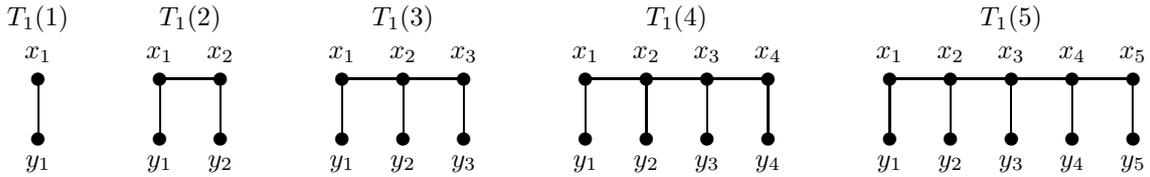

For a positive integer $k$, let $T_1(k)$ be the tree illustrated in Figure \ref{fig1}, that is,
$T_1(k)$ has vertex set $\{ x_1,\ldots,x_k\}\cup \{ y_1,\ldots,y_k\}$,
contains the path $x_1\ldots x_k$,
and $x_i$ is the only neighbor of $y_i$ for $i\in [k]$.

\begin{figure}[H]
\begin{center}
{\footnotesize
$\mbox{}$\hfill
%TeXCAD Picture [2.pic]. Options:
%\grade{\on}
%\emlines{\off}
%\epic{\off}
%\beziermacro{\on}
%\reduce{\on}
%\snapping{\on}
%\pvinsert{% Your \input, \def, etc. here}
%\quality{8.000}
%\graddiff{0.005}
%\snapasp{1}
%\zoom{11.3137}
\unitlength 0.8mm % = 2.845pt
\linethickness{0.4pt}
\ifx\plotpoint\undefined\newsavebox{\plotpoint}\fi % GNUPLOT compatibility
\begin{picture}(36,25)(0,0)
\put(25,5){\circle*{2}}
\put(15,5){\circle*{2}}
\put(35,5){\circle*{2}}
\put(25,15){\circle*{2}}
\put(15,15){\circle*{2}}
\put(35,15){\circle*{2}}
\put(25,5){\line(0,1){10}}
\put(15,5){\line(0,1){10}}
\put(35,5){\line(0,1){10}}
\put(15,15){\line(1,0){10}}
\put(25,15){\line(1,0){10}}
\put(25,19){\makebox(0,0)[cc]{$x_2$}}
\put(15,19){\makebox(0,0)[cc]{$x_1$}}
\put(35,19){\makebox(0,0)[cc]{$x_3$}}
\put(25,1){\makebox(0,0)[cc]{$y_2$}}
\put(15,1){\makebox(0,0)[cc]{$y_1$}}
\put(35,1){\makebox(0,0)[cc]{$y_3$}}
\put(20,25){\makebox(0,0)[cc]{$T_2(3)$}}
\put(5,15){\circle*{2}}
\put(15,15){\line(-1,0){10}}
\put(5,19){\makebox(0,0)[cc]{$a$}}
\end{picture}
\hfill
%TeXCAD Picture [3.pic]. Options:
%\grade{\on}
%\emlines{\off}
%\epic{\off}
%\beziermacro{\on}
%\reduce{\on}
%\snapping{\on}
%\pvinsert{% Your \input, \def, etc. here}
%\quality{8.000}
%\graddiff{0.005}
%\snapasp{1}
%\zoom{11.3137}
%\unitlength 1mm % = 2.845pt
\linethickness{0.4pt}
\ifx\plotpoint\undefined\newsavebox{\plotpoint}\fi % GNUPLOT compatibility
\begin{picture}(66,25)(0,0)
\put(55,5){\circle*{2}}
\put(45,5){\circle*{2}}
\put(65,5){\circle*{2}}
\put(55,15){\circle*{2}}
\put(25,15){\circle*{2}}
\put(45,15){\circle*{2}}
\put(15,15){\circle*{2}}
\put(65,15){\circle*{2}}
\put(55,5){\line(0,1){10}}
\put(45,5){\line(0,1){10}}
\put(65,5){\line(0,1){10}}
\put(45,15){\line(1,0){10}}
\put(15,15){\line(1,0){10}}
\put(55,15){\line(1,0){10}}
\put(25,15){\line(1,0){10}}
\put(55,19){\makebox(0,0)[cc]{$x_2$}}
\put(45,19){\makebox(0,0)[cc]{$x_1$}}
\put(65,19){\makebox(0,0)[cc]{$x_3$}}
\put(55,1){\makebox(0,0)[cc]{$y_2$}}
\put(45,1){\makebox(0,0)[cc]{$y_1$}}
\put(65,1){\makebox(0,0)[cc]{$y_3$}}
\put(35,25){\makebox(0,0)[cc]{$T_3(3)$}}
\put(35,15){\circle*{2}}
\put(5,15){\circle*{2}}
\put(45,15){\line(-1,0){10}}
\put(15,15){\line(-1,0){10}}
\put(35,19){\makebox(0,0)[cc]{$a$}}
\put(25,19){\makebox(0,0)[cc]{$b$}}
\put(15,19){\makebox(0,0)[cc]{$c$}}
\put(5,19){\makebox(0,0)[cc]{$d$}}
\end{picture}
\hfill
%TeXCAD Picture [4.pic]. Options:
%\grade{\on}
%\emlines{\off}
%\epic{\off}
%\beziermacro{\on}
%\reduce{\on}
%\snapping{\on}
%\pvinsert{% Your \input, \def, etc. here}
%\quality{8.000}
%\graddiff{0.005}
%\snapasp{1}
%\zoom{11.3137}
%\unitlength 1mm % = 2.845pt
\linethickness{0.4pt}
\ifx\plotpoint\undefined\newsavebox{\plotpoint}\fi % GNUPLOT compatibility
\begin{picture}(46,30)(0,0)
\put(35,5){\circle*{2}}
\put(25,5){\circle*{2}}
\put(45,5){\circle*{2}}
\put(35,15){\circle*{2}}
\put(5,25){\circle*{2}}
\put(25,15){\circle*{2}}
\put(45,15){\circle*{2}}
\put(35,5){\line(0,1){10}}
\put(25,5){\line(0,1){10}}
\put(45,5){\line(0,1){10}}
\put(25,15){\line(1,0){10}}
\put(35,15){\line(1,0){10}}
\put(5,25){\line(1,0){10}}
\put(35,19){\makebox(0,0)[cc]{$x_2$}}
\put(21,15){\makebox(0,0)[cc]{$x_1$}}
\put(45,19){\makebox(0,0)[cc]{$x_3$}}
\put(35,1){\makebox(0,0)[cc]{$y_2$}}
\put(25,1){\makebox(0,0)[cc]{$y_1$}}
\put(45,1){\makebox(0,0)[cc]{$y_3$}}
\put(25,30){\makebox(0,0)[cc]{$T_4(3)$}}
\put(15,25){\circle*{2}}
\put(15,29){\makebox(0,0)[cc]{$a$}}
\put(5,29){\makebox(0,0)[cc]{$b$}}
\put(15,25){\line(2,-1){20}}
\end{picture}
\hfill$\mbox{}$\\[10mm]
%TeXCAD Picture [5.pic]. Options:
%\grade{\on}
%\emlines{\off}
%\epic{\off}
%\beziermacro{\on}
%\reduce{\on}
%\snapping{\on}
%\pvinsert{% Your \input, \def, etc. here}
%\quality{8.000}
%\graddiff{0.005}
%\snapasp{1}
%\zoom{11.3137}
%\unitlength 1mm % = 2.845pt
\linethickness{0.4pt}
\ifx\plotpoint\undefined\newsavebox{\plotpoint}\fi % GNUPLOT compatibility
\begin{picture}(106,25)(0,0)
\put(55,5){\circle*{2}}
\put(45,5){\circle*{2}}
\put(65,5){\circle*{2}}
\put(55,15){\circle*{2}}
\put(25,15){\circle*{2}}
\put(95,15){\circle*{2}}
\put(45,15){\circle*{2}}
\put(15,15){\circle*{2}}
\put(85,15){\circle*{2}}
\put(65,15){\circle*{2}}
\put(55,5){\line(0,1){10}}
\put(45,5){\line(0,1){10}}
\put(65,5){\line(0,1){10}}
\put(45,15){\line(1,0){10}}
\put(15,15){\line(1,0){10}}
\put(85,15){\line(1,0){10}}
\put(55,15){\line(1,0){10}}
\put(25,15){\line(1,0){10}}
\put(95,15){\line(1,0){10}}
\put(55,19){\makebox(0,0)[cc]{$x_2$}}
\put(45,19){\makebox(0,0)[cc]{$x_1$}}
\put(65,19){\makebox(0,0)[cc]{$x_3$}}
\put(55,1){\makebox(0,0)[cc]{$y_2$}}
\put(45,1){\makebox(0,0)[cc]{$y_1$}}
\put(65,1){\makebox(0,0)[cc]{$y_3$}}
\put(55,25){\makebox(0,0)[cc]{$T_5(3)$}}
\put(35,15){\circle*{2}}
\put(105,15){\circle*{2}}
\put(5,15){\circle*{2}}
\put(75,15){\circle*{2}}
\put(45,15){\line(-1,0){10}}
\put(15,15){\line(-1,0){10}}
\put(85,15){\line(-1,0){10}}
\put(35,19){\makebox(0,0)[cc]{$a$}}
\put(105,19){\makebox(0,0)[cc]{$d'$}}
\put(25,19){\makebox(0,0)[cc]{$b$}}
\put(95,19){\makebox(0,0)[cc]{$c'$}}
\put(15,19){\makebox(0,0)[cc]{$c$}}
\put(85,19){\makebox(0,0)[cc]{$b'$}}
\put(5,19){\makebox(0,0)[cc]{$d$}}
\put(75,19){\makebox(0,0)[cc]{$a'$}}
\put(65,15){\line(1,0){10}}
\end{picture}
}
\end{center}
\caption{The trees $T_k(3)$ for $k\in \{ 2,3,4,5\}$.}\label{fig2}
\end{figure}
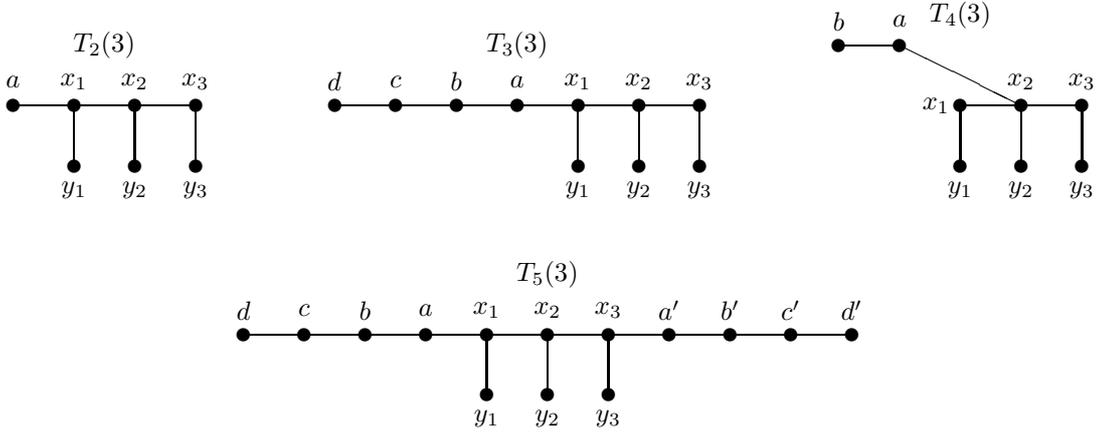

\noindent Let $T_2(k)$ arise from $T_1(k)$ by adding a vertex $a$ and the edge $x_1a$.
Let $T_3(k)$ arise from $T_1(k)$ by adding the vertices $a$, $b$, $c$, and $d$, and the edges $x_1a$, $ab$, $bc$, and $cd$.
For $k\geq 3$, let $T_4(k)$ arise from $T_1(k)$ by adding the vertices $a$ and $b$, and the edges $x_2a$ and $ab$.
Finally, let $T_5(k)$ arise from $T_1(k)$ by adding the vertices
$a$, $b$, $c$, $d$, $a'$, $b'$, $c'$, and $d'$,
and the edges
$x_1a$, $ab$, $bc$, $cd$, $x_ka'$, $a'b'$, $b'c'$, $c'd'$.
See Figure \ref{fig2} for an illustration.

\noindent Let
$${\cal T}=\{ P_1,P_8\}\cup\bigcup_{k\in \mathbb{N}}\big\{ T_1(k),T_2(k),T_3(k),T_5(k)\big\}
\cup\bigcup_{k\geq 3}\big\{ T_4(k)\big\}.$$
Note that ${\cal T}$ contains the paths 
$P_1$,
$P_2=T_1(1)$, 
$P_3=T_2(1)$, 
$P_4=T_1(2)$, 
$P_6=T_3(1)$, and
$P_8$.

\begin{lemma}\label{lemma2}
Every tree $T\in {\cal T}$ satisfies $\alpha_e(T)=\alpha(T)$.
Furthermore, if $S$ is a maximum exponential independent set of $T$, then
\begin{enumerate}[(i)]
\item $S\in \big\{ 
\{ y_1,\ldots,y_k\},
\{ x_1\}\cup \{ y_2,\ldots,y_k\},
\{ x_k\}\cup \{ y_1,\ldots,y_{k-1}\}
\big\}$,
\item $S=\{ a,y_1,\ldots,y_k\}$ if $T=T_2(k)$,
\item $S=\{ b,d,y_1,\ldots,y_k\}$ if $T=T_3(k)$ with $k\geq 2$,
\item $S=\{ b,y_1,\ldots,y_k\}$ if $T=T_4(k)$ with $k\geq 3$,
\item $S=\{ b,d,b',d',y_1,\ldots,y_k\}$ if $T=T_5(k)$, 
\item $S\in \big\{ \{ y_1,a,d\} ,\{ y_1,b,d\}\big\}$ if $T=P_6=y_1x_1abcd$, and 
\item $S=\{ b,d,b',d'\}$ if $T=P_8=dcbaa'b'c'd'$.
\end{enumerate}
\end{lemma}
{\it Proof:} Let $T\in {\cal T}$.
It is easy to see that $\alpha_e(P_n)=\alpha(P_n)$ for $n\in \{ 1,2,3,4,6,8\}$.
Furthermore, one easily checks that $P_6$ has only two distinct exponential independent sets of order $3$,
and that $P_8$ has a unique exponential independent set of order $4$,
which implies (vi) and (vii).
Now, we may assume that $T\not\in \{ P_1,P_6,P_8\}$.

First, we assume that $T=T_1(k)$ for some positive integer $k$.
Clearly, the set $\{ y_1,\ldots,y_k\}$ is a maximum independent set, which implies $\alpha(T)=k$.
Since this set is also exponential independent, we obtain $\alpha_e(T)=\alpha(T)=k$.
Now, let $S$ be a maximum exponential independent set of $T$.
For $k\in [2]$, it follows easily that $S$ is as stated in (i). 
Now, let $k\geq 3$.
Since $S$ contains at most one of the two vertices $x_i$ and $y_i$ for each $i\in [k]$,
the set $S$ necessarily intersects each of the sets $\{ x_i,y_i\}$ for $i\in [k]$ in exactly one vertex.
If $x_i\in S$ for some $i\in \{ 2,\ldots,k-1\}$,
this implies that $y_{i-1},y_{i+1}\in S$, 
which yields the contradiction that $w_{(T,S\setminus \{ x_i\})}(x_i)\geq \frac{1}{2}+\frac{1}{2}=1$.
Hence, $\{ y_i:2\leq i\leq k-1\}\subseteq S$.
If $x_1,x_k\in S$, 
then $w_{(T,S\setminus \{ x_1\})}(x_1)=\sum\limits_{i=2}^{k-1}\left(\frac{1}{2}\right)^{i-1}+\frac{1}{2^{k-2}}=1$,
which is a contradiction.
Hence, the set $S$ is stated as in (i).

Next, we assume that $T=T_2(k)$ for some positive integer $k$.
Again, the set of leaves is a maximum independent set of $T$,
which is also exponential independent, 
and, hence, $\alpha_e(T)=\alpha(T)=k+1$.
Now, let $S$ be a maximum exponential independent set of $T$.
If $a\not\in S$, then $S$ is an exponential independent set of $T-a=T_1(k)$,
which contradicts $\alpha_e(T_1(k))=k$. Hence, $a\in S$, which implies $x_1\not\in S$.
For $k\in [2]$, it follows easily that $S$ is as stated in (ii).
Now, let $k\geq 3$.
Since $S\setminus \{ a\}$ is a maximum exponential independent set of $T-a=T_1(k)$, we obtain, by (i), that $\{ y_1,\ldots,y_{k-1}\}\subseteq S$.
If $x_k\in S$, then $w_{(T,S\setminus \{ a\})}(a)\geq 1$ follows similarly as above,
which is a contradiction.
Hence, the set $S$ is as stated in (ii).

Next, we assume that $T=T_3(k)$ for some positive integer $k$.
Since $T\not=P_6$, we have $k\geq 2$.
As before, it follows easily that the set specified in (iii) is a maximum exponential independent set of $T$,
and, hence, $\alpha_e(T)=\alpha(T)=k+2$.
Now, let $S$ be a maximum exponential independent set of $T$.
Necessarily, $|S\setminus \{ a,b,c,d\}|=k$ and $|S\cap \{ a,b,c,d\}|=2$,
which implies that $S$ contains either $a$ or $b$.
If $S$ contains $a$, then $S\setminus \{ b,c,d\}$ is a maximum exponential independent set of $T_2(k)$,
which, by (ii), implies $S\setminus \{ b,c,d\}=\{ a,y_1,\ldots,y_k\}$.
Now, we obtain the contradiction,
$w_{(T,S\setminus \{ a\})}(a)
=
w_{(T,S\setminus \{ a,b,c,d\})}(a)
+
w_{(T,S\cap \{ b,c,d\})}(a)
\geq 
\frac{1}{2}+\frac{1}{4}+\frac{1}{4}=1$.
Hence, $b\in S$, which implies $S\cap \{ a,b,c,d\}=\{ b,d\}$, 
and $x_1\not\in S$.
Since $S\setminus \{ a,b,c,d\}$ is a maximum exponential independent set of $T-\{ a,b,c,d\}=T_1(k)$,
we obtain, by (i), that $\{ y_1,\ldots,y_{k-1}\}\subseteq S$.
If $x_k\in S$, then $w_{(T,S\setminus \{ a\})}(a)=1$, which is a contradiction.
Hence, the set $S$ is as stated in (iii).

Finally, if 
either 
$T=T_4(k)$ for some integer $k$ with $k\geq 3$ 
or
$T=T_5(k)$ for some positive integer $k$,
very similar arguments as above imply that $\alpha_e(T)=\alpha(T)$, 
and that every maximum exponential independent set is as specified in (iv) and (v).
$\Box$ 

\begin{theorem}\label{theorem6}
If $T$ is a tree, then $\alpha_e(T)=\alpha(T)$ if and only if $T\in {\cal T}$.
\end{theorem}
{\it Proof:} In view of Lemma \ref{lemma2}, 
it remains to show that every tree $T$ with $\alpha_e(T)=\alpha(T)$ belongs to ${\cal T}$.
Therefore, suppose that $T$ is a tree of minimum order such that $\alpha_e(T)=\alpha(T)$ but $T\not\in{\cal T}$.
Let $S$ be a maximum exponential independent set of $T$.
Since 
$P_1,P_2,P_3\in {\cal T}$,
and 
$\alpha_e(K_{1,n-1})=2<n-1=\alpha(K_{1,n-1})$ for $n\geq 3$,
we may assume that $T$ has diameter at least $3$.
Therefore, if $w$ is an endvertex of a longest path in $T$,
then the unique neighbor $v$ of $w$ has exactly one neighbor $u$ that is not an endvertex of $T$.
Let $T'$ be the component of $T-v$ that contains $u$.
Note that $T'$ is not $P_1$.
We consider different cases.

\medskip

\noindent {\bf Case 1} $d_T(v)\geq 4$.

\medskip

\noindent Clearly, $\alpha(T)\geq \alpha(T')+(d_T(v)-1)\geq \alpha(T')+3$.
Since $S$ contains at most $2$ vertices from $N_T[v]\setminus \{ u\}$,
and $S\cap V(T')$ is an exponential independent set of $T'$, 
we obtain $\alpha_e(T)\leq \alpha_e(T')+2$,
which yields the contradiction
$\alpha(T)\geq \alpha(T')+3\geq \alpha_e(T')+3>\alpha_e(T)$,
which completes the proof in this case.

\medskip

\noindent {\bf Case 2} $d_T(v)=3$.

\medskip

\noindent Let $N_T(v)=\{ u,w,w'\}$.
As before, we obtain that $\alpha(T)\geq \alpha(T')+2$ and $\alpha_e(T)\leq \alpha_e(T')+2$,
which implies
$\alpha(T)\geq \alpha(T')+2\geq \alpha_e(T')+2\geq \alpha_e(T)=\alpha(T)$.
Since equality holds throughout this inequality chain,
we have 
$\alpha(T')=\alpha_e(T')$
and 
$\alpha_e(T)=\alpha_e(T')+2$.
By the choice of $T$, the condition $\alpha(T')=\alpha_e(T')$ implies that $T'\in {\cal T}$.
Furthermore, $\alpha_e(T)=\alpha_e(T')+2$ implies that $S\cap \{ v,w,w'\}=\{ w,w'\}$, and that 
$S'=S\setminus \{ w,w'\}$ is a maximum exponential independent set of $T'$.
Since $w_{(T,\{ w,w'\})}(u)=1$, we obtain $u\not\in S'$.

First, we assume that $T'=P_8=dcbaa'b'c'd'$.
By Lemma \ref{lemma2}, we have $S'=\{ b,d,b',d'\}$.
By symmetry, we may assume that $u\in \{ a,c\}$.
In both cases $w_{(T,\{ d,w,w'\})}(b)=1$,
which is a contradiction.

Next, we assume that $T'=T_1(k)$.
Since $u\not\in S'$,
we have $u=x_i$ for some $i\in [k]$.
If $i\in \{ 1,k\}$, then $T=T_2(k+1)\in {\cal T}$, 
which is a contradiction. 
Hence, $2\leq i\leq k-1$.
Using Lemma \ref{lemma2}, 
we obtain the contradiction
$w_{(T,S\setminus \{ y_i\})}(y_i)\geq 4\cdot \frac{1}{4}$.

Next, we assume that $T'=T_2(k)$.
By Lemma \ref{lemma2},
we have $S'=\{ a,y_1,\ldots,y_k\}$.
Since $u\not\in S'$,
we have $u=x_i$ for some $i\in [k]$,
which implies the contradiction
$w_{(T,S\setminus\{ y_i\})}(y_i)\geq 
w_{(T,\{ w,w'\})}(y_i)
+
w_{(T,\{ a,y_1,\ldots,y_{i-1}\})}(y_i)
=
2\cdot\frac{1}{4}
+\left(\frac{1}{4}+\cdots+\frac{1}{2^i}+\frac{1}{2^i}\right)=1$.

Next, we assume that $T'=T_3(k)$.
By Lemma \ref{lemma2},
we have $S'=\{ b,d,y_1,\ldots,y_k\}$.
Since $u\not\in S'$,
we have $u\in \{ a,c\}$ or $u=x_i$ for some $i\in [k]$.
In the former case,
$w_{(T,\{ d,w,w'\})}(b)=1$, 
and in the latter case,
$w_{(T,S\setminus\{ b\})}(b)\geq 
w_{(T,\{ d,w,w'\})}(b)
+
w_{(T,\{ y_1,\ldots,y_i\})}(b)
=
\frac{1}{2}+2\cdot\frac{1}{2^{i+2}}
+\left(\frac{1}{4}+\cdots+\frac{1}{2^{i+1}}\right)=1$.

Next, we assume that $T'=T_4(k)$ for some $k\geq 3$.
By Lemma \ref{lemma2},
we have $S'=\{ b,y_1,\ldots,y_k\}$.
Since $u\not\in S'$,
we have $u\in \{ a,x_1\}$ or $u=x_i$ for some $i\in \{ 2,\ldots,k\}$.
If $u=a$, 
then 
$w_{(T,S\setminus \{ b\})}(b)\geq 
w_{(T,\{ w,w',y_1,y_2,y_3\})}(b)
=3\cdot \frac{1}{4}+2\cdot \frac{1}{8}=1$.
Similarly, if $u=x_1$,
then $w_{(T,S\setminus \{ y_1\})}(y_1)\geq 1$.
Finally, if $u=x_i$ for some $i\geq 2$, then
$w_{(T,S\setminus \{ y_i\})}(y_i)\geq 
w_{(T,\{ w,w',b\})}(y_i)
+
w_{(T,\{ y_1,\ldots,y_{i-1}\})}(y_i)
=
2\cdot\frac{1}{4}+\frac{1}{2^i}
+\left(\frac{1}{4}+\cdots+\frac{1}{2^i}\right)=1$.

Finally, we assume that $T'=T_5(k)$.
By Lemma \ref{lemma2},
we have $S'=\{ b,d,b',d',y_1,\ldots,y_k\}$.
Since $u\not\in S'$,
we have $u\in \{ a,c,a',c'\}$ or $u=x_i$ for some $i\in \{ 1,\ldots,k\}$.
If 
$u\in \{ a,c\}$, then
$w_{(T,S\setminus \{ b\})}(b)\geq 
w_{(T,\{ d,w,w'\})}(b)=1$,
if
$u\in \{ a',c'\}$, then
$w_{(T,S\setminus \{ b'\})}(b')\geq 1$, and
if $u=x_i$, then
$w_{(T,S\setminus \{ b\})}(b)\geq 
w_{(T,\{ d,w,w'\})}(b)
+
w_{(T,\{ y_1,\ldots,y_i\})}(b)=
\frac{1}{2}+2\cdot\frac{1}{2^{i+2}}
+\left(\frac{1}{4}+\cdots+\frac{1}{2^{i+1}}\right)=1$,
which completes the proof in this case.

\medskip

\noindent {\bf Case 3} $d_T(v)=2$.

\medskip

\noindent Let $N_T(v)=\{ u,w\}$.
As before, we obtain that $\alpha(T)\geq \alpha(T')+1$ and $\alpha_e(T)\leq \alpha_e(T')+1$,
which implies
$\alpha(T)\geq \alpha(T')+1\geq \alpha_e(T')+1\geq \alpha_e(T)=\alpha(T)$.
Again, equality holds throughout this inequality chain,
and we obtain that
$\alpha(T')=\alpha_e(T')$,
$\alpha_e(T)=\alpha_e(T')+1$,
$T'\in {\cal T}$, and 
$S'=S\setminus \{ v,w\}$ is a maximum exponential independent set of $T'$.
Clearly, we may assume that $S\cap \{ v,w\}=\{ w\}$.

First, we assume that $T'=P_8=dcbaa'b'c'd'$.
By Lemma \ref{lemma2}, we have $S'=\{ b,d,b',d'\}$.
In each case, either $w_{(T,S\setminus \{ b\})}(b)\geq 1$
or $w_{(T,S\setminus \{ b'\})}(b')\geq 1$.

Next, we assume that $T'=T_1(k)$.
If $u\in \{ x_1,x_k\}$, then $T=T_1(k+1)\in {\cal T}$.
If $u\in \{ y_1,y_k\}$, then 
either $k=1$ and $T=P_4\in {\cal T}$,
or $k\geq 2$ and $T=T_3(k-1)\in {\cal T}$.
If $u\in \{ x_2,x_{k-1}\}$, then $T=T_4(k)\in {\cal T}$.
If $u=y_i$ for some $i\in \{ 2,\ldots,k-1\}$,
then, by Lemma \ref{lemma2}, 
$y_i\in S'$ and 
$w_{(T,S\setminus \{ y_i\})}(y_i)
\geq w_{(T,\{ w,y_{i-1},y_{i+1}\})}(y_i)=1$.
Finally,
if $u=x_i$ for some $i\in \{ 3,\ldots,k-2\}$,
then, by Lemma \ref{lemma2}, 
$y_i\in S'$ and 
$w_{(T,S\setminus \{ y_i\})}(y_i)
\geq w_{(T,\{ w,y_{i-1},y_{i+1},y_{i-2},y_{i+2}\})}(y_i)=1$.

Next, we assume that $T'=T_2(k)$.
If $u\in \{ a,y_1,\ldots,y_k\}$, then, by Lemma \ref{lemma2},
$w_{(T,S\setminus \{ u\})}(u)=
\frac{1}{2}+w_{(T,S'\setminus \{ u\})}(u)
\geq
\frac{1}{2}+\frac{1}{2}=1$.
If $u=x_k$, then $T=T_2(k+1)\in {\cal T}$.
Finally,
if $u=x_i$ for some $i\in [k-1]$, then,
by Lemma \ref{lemma2}, 
$y_i\in S'$ and
$w_{(T,S\setminus \{ y_i\})}(y_i)
\geq 
w_{(T,\{ w,y_{i+1}\})}(y_i)
+
w_{(T,\{ b,y_1,\ldots,y_{i-1}\})}(y_i)
=2\cdot \frac{1}{4}
+\left(\frac{1}{4}+\cdots+\frac{1}{2^i}+\frac{1}{2^i}\right)=1$.

Next, we assume that $T'=T_3(k)$.
If $u\in \{ a,b,c,d\}$, then, by Lemma \ref{lemma2},
$w_{(T,S\setminus \{ b\})}(b)\geq 
w_{(T,\{ w\})}(b)
+
w_{(T,\{ d,y_1\})}(b)\geq 1$.
If $u=x_k$, then $T=T_3(k+1)\in {\cal T}$.
If $u=y_k$, then 
either $k=1$ and $T=P_8\in {\cal T}$,
or $k\geq 2$ and $T=T_5(k-1)\in {\cal T}$.
If $u=y_i$ for some $i\in [k-1]$,
then $y_i\in S'$ and
$w_{(T,S\setminus \{ y_i\})}(y_i)\geq \frac{1}{2}+2\cdot\frac{1}{4}=1$.
Finally,
if $u=x_i$ for some $i\in [k-1]$,
then 
$w_{(T,S\setminus \{ b\})}(b)\geq 
w_{(T,\{ d,w\})}(b)
+
w_{(T,\{ y_1,\ldots,y_{i+1}\})}(b)
=\frac{1}{2}+\frac{1}{2^{i+2}}
+\left(\frac{1}{4}+\cdots+\frac{1}{2^{i+2}}\right)=1$.

Next, we assume that $T'=T_4(k)$ for some $k\geq 3$.
If $u\in \{ b,y_1,\ldots,y_k\}$, then, by Lemma \ref{lemma2},
$u\in S'$ and 
$w_{(T,S\setminus \{ u\})}(u)=
\frac{1}{2}+w_{(T,S'\setminus \{ u\})}(u)
\geq
\frac{1}{2}+\frac{1}{2}=1$.
If $u=x_1$ and $k=3$, then $T=T_4(4)\in {\cal T}$.
If $u=x_1$ and $k\geq 4$, then 
$w_{(T,S\setminus \{ y_2\})}(y_2)\geq
w_{(T,\{ b,w,y_1,y_3,y_4\})}(y_2)
=3\cdot\frac{1}{4}+2\cdot\frac{1}{8}=1$.
If $u=a$, we obtain similar contradictions.
If $u=x_k$, then $T=T_4(k+1)\in {\cal T}$.
Finally, if $u=x_i$ for some $i\in \{ 2,\ldots,k-1\}$, then
$y_i\in S'$ and
$w_{(T,S\setminus \{ y_i\})}(y_i)\geq
w_{(T,\{ y_{i+1},w\})}(y_i)
+w_{(T,\{ b,y_1,\ldots,y_{i-1},w\})}(y_i)=
2\cdot \frac{1}{4}
+\left(\frac{1}{4}+\cdots+\frac{1}{2^i}+\frac{1}{2^i}\right)=1$.

Finally, we assume that $T'=T_5(k)$.
If $u\in \{ b,d,b',d',y_1,\ldots,y_k\}$, then, by Lemma \ref{lemma2},
$u\in S'$ and 
$w_{(T,S\setminus \{ u\})}(u)=
\frac{1}{2}+w_{(T,S'\setminus \{ u\})}(u)
\geq
\frac{1}{2}+\frac{1}{2}=1$.
If $u\in \{ a,c\}$, then 
$w_{(T,S\setminus \{ b\})}(b)\geq w_{(T,\{ d,w,y_1\})}(b)=1$.
If $u\in \{ a',c'\}$, 
we obtain a similar contradiction.
Finally, if $u=x_i$ for some $i\in [k]$, then 
$w_{(T,S\setminus \{ b\})}(b)\geq
w_{(T,\{ w,d\})}(b)
+w_{(T,\{ y_j:1\leq j\leq \min\{ i+1,k\}\}\cup \{ b'\})}(b)=
\frac{1}{2}+\frac{1}{2^{i+2}}
+\left(\frac{1}{4}+\cdots+\frac{1}{2^{i+2}}\right)=1$,
which completes the proof.
$\Box$

\section{Conclusion}

Our results motivate several open problems.
It seems interesting to characterize all extremal graphs for Theorem \ref{theorem2}.
In view of Theorem \ref{theorem2b}, one can study upper bounds for graphs of larger minimum degree.
As stated before Theorem \ref{theorem4}, we conjecture that $\alpha_e(G)$ grows faster than $\log_2(n(G))$ for cubic graphs.
Can the graphs $G$ with $\alpha_e(G)=\alpha(G)$ be recognized efficiently?
Are there hardness results concerning $\alpha_e(G)$, and efficient algorithms for restricted graphs classes such as trees?
% porous version
% algorithm for subcubic trees

\end{document}